\documentclass[12pt]{article}
\usepackage{amsmath, amsthm, amssymb}
\usepackage{amsfonts}
\usepackage{epsfig}
\usepackage{graphics}
\usepackage{subfigure}
\usepackage{graphicx}
\usepackage{color}
\usepackage{epstopdf} 
\usepackage{float} 
\usepackage{mathrsfs} 
\usepackage{relsize} 
\newcommand{\bq}{\begin{equation}}
\newcommand{\eq}{\end{equation}}
\newcommand{\bqr}{\begin{eqnarray}}
\newcommand{\eqr}{\end{eqnarray}}
\newcommand{\bqrn}{\begin{eqnarray*}}
\newcommand{\eqrn}{\end{eqnarray*}}
\begin{document}

\noindent {\large \bf Mathematics of Parking}\\
{\bf Varying Parking Rate}\\[2mm]
Pavel B. Dubovski, Michael Tamarov\\
Published in {\it Journal of Statistical Physics} v. {\bf 182}, no. 22 (2021)\\
DOI: 10.1007/s10955-020-02678-x\\

\begin{abstract}
In the classical parking problem, unit intervals ("car lengths") are placed uniformly at random without overlapping. The process terminates at saturation, i.e. until no more unit intervals can be stowed. In this paper, we present a generalization of this problem in which the unit intervals are placed with an exponential distribution with rate parameter $\lambda$. We show that the mathematical expectation of the number of intervals present at saturation satisfies a certain integral equation. Using Laplace transforms and Tauberian theorems, we investigate the asymptotic behavior of this function and describe a way to compute the corresponding limits for large $\lambda$. Then, we derive another integral equation for the derivative of this function and use it to compute the above limits for small $\lambda$ with the help of some asymptotic results for integral equations. We also show that the corresponding limits converge to the uniform case as $\lambda$ vanishes, yielding the well-known Renyi constant. Finally, we reveal the asymptotic behavior of the variance of the intervals at saturation.

\noindent{\bf keywords:} asymptotic behavior; integral equations; parking problem; percolation; R\'{e}nyi constant; saturation; self-similar random variable; Tauberian theorems

\noindent MSC 60K35, 82B43,  45D99
\end{abstract}

\section{Introduction}
\label{intro}
In his seminal paper, Alfred R\'{e}nyi \cite{Renyi} presents and solves the following parking problem. For $x\ge0$, place an open unit interval at random on the interval $(0,x)$ by choosing the left endpoint of the unit interval uniformly at random from the interval $(0,x-1)$. Then, randomly place a second open unit interval on $(0,x)$, independent of the first. If it intersects the first unit interval, discard it and repeat until the two intervals are disjoint. In general, if $k$ disjoint unit intervals have already been placed, choose the next interval uniformly at random on the interval $(0,x)$, discarding it if it intersects any of the $k$ previously chosen intervals. The process is repeated until saturation, i.e. until no more unit intervals can be placed without intersecting any of the ones already placed. Letting the random variable $v_x$ denote the number of unit intervals present at saturation and $M(x)=E[v_x]$, the expectation of $v_x$, R\'{e}nyi  shows that $M(x)$ satisfies 
\begin{equation}
M(x+1)=\frac{2}{x}\int_0^xM(t)\,dt+1,\label{eq:1}\end{equation}
$$M(x)=0\quad\mbox{for }0\le x\le1.$$

\noindent Then, using Laplace transforms and a Tauberian theorem, he proves that 
$$\lim_{x\to\infty}\frac{M(x)}{x}=C,$$
where 
$$C=\int_0^\infty\exp\left(-2\int_0^t\frac{1-e^{-u}}{u}\right)\,dt\approx0.748$$
and moreover that \begin{equation}M(x)=Cx-(1-C)+O\left(x^{-n}\right)\label{eq:2}\end{equation}
for each $n=1,2,\ldots$. Using some asymptotic results for integral equations, Dvoretsky and Robbins \cite{Dvoretsky} later improve the bound on the error term to 
\begin{equation}O\left(\left(\frac{2e}{x}\right)^{x-3/2}\right).\label{Dvoretsky1}\end{equation}
R\'{e}nyi also shows that Var$(v_x)=O(x)$, a result that Blaisdell and Solomon \cite{Blaisdell} later sharpen, showing that 
$$\lim_{x\to\infty}\frac{\mbox{Var}(v_x)}{x}=D,$$
where $D\approx0.035672$, and which Dvoretsky and Robbins \cite{Dvoretsky} subsequently improve to
\begin{equation}
\mbox{Var}(v_x)=Dx+D+O\left(\left(\frac{4e}{x}\right)^{x-4}\right).\label{Dvoretsky2}
\end{equation}

Since R\'{e}nyi published his paper, many versions of this basic problem have been considered. For example, Mullooly \cite{Mullooly} presents a generalization in which both the position and length of the smaller intervals are independent, uniform random variables (with the length of the interval bounded below by some positive constant). Similarly, Krapivsky \cite{Krapivsky1} considers the case where intervals are placed on an infinite line, where the parking distribution varies as a power of the length of the interval. In a discrete version of the problem, Ziff \cite{Ziff1} studies the distribution of gaps in a random sequential adsorption process of dimers on a one-dimensional lattice, and in a slightly different model, Mansfield \cite{Mansfield} investigates the random attachment of probe spheres of radius $r_1$ on a target sphere of radius $r_2$, where the probe spheres are placed on the surface of the target sphere uniformly at random without overlapping. He then considers the case where the probe spheres diffuse into place and shows that this latter model produces a higher packing density than the former. 

In this paper, we consider another generalization in which the unit intervals are not chosen uniformly at random, but with a truncated exponential distribution with rate parameter $\lambda>0$ and support $[0,x-1]$. In other words, at each step of the selection process, the left endpoint of the next unit interval is chosen with a probability distribution whose density function is given by 
\begin{equation}
f(t;\lambda,x)=\begin{cases}\displaystyle\frac{\lambda e^{-\lambda t}}{1-e^{-\lambda (x-1)}},&\quad0<t< x-1\vspace{1mm}\\0,&\quad\mbox{otherwise.}\end{cases}\label{200}
\end{equation}
This corresponds to cars with a preference to park closer to some popular destination located at one end of a street. In fact, this is related to a paper by Krapivsky and Redner \cite{Krapivsky2} in which they examine various parking strategies in this situation, with cars entering one side of the street at some rate $\lambda$ and leaving at a unit rate.

R\'{e}nyi's brilliant approach in the original parking problem is not easily generalizable for arbitrary probability distributions, but it can be adapted to this case because the (truncated) exponential distribution $f(t;\lambda,x)$ is self-similar in the following sense.\\

\noindent{\bf Definition.} Let $X$ be a random variable with probability density function 
$f_X(t)$ and let $Y_{x_1,x_2}$ be the conditional random variable, 
$Y_{x_1,x_2}=X|\left(x_1\le X\le x_2\right)$. We call the random variable $X$ 
{\it self-similar on the interval $[a,b]$} if for all $a\le x_1 < x_2\le b$, 
\[
f_Y(t)=\begin{cases}C\cdot f_X(t-x_1),&\quad x_1\le t\le x_2\\
0,&\quad\mbox{otherwise,}\end{cases}
\]
where $C=\left(\displaystyle\int_{x_1}^{x_2}f_X(t)\,dt\right)^{-1}$ is the normalizing constant.
Thus, the random variable $X$ is self-similar if the conditional random variable $Y$ can be obtained from $X$ with a shift, normalization, and restriction to the appropriate interval.

The truncated exponential random variable $X$ defined by \eqref{200} is self-similar on the interval $[0,x-1]$. In fact, for any $0\le x_1<x_2\le x-1$,
\begin{eqnarray*}
f_Y=f_{x_1,x_2}(t;\lambda,x)&=&\frac{f_X(t)}{\displaystyle\int_{x_1}^{x_2}f_X(t)\,dt}
=\frac{\displaystyle\left(\frac{\lambda e^{-\lambda t}}{1-e^{-\lambda(x-1)}}\right)}
{\displaystyle\left(\frac{e^{-\lambda x_1}-e^{-\lambda x_2}}{1-e^{-\lambda(x-1)}}\right)}\\
&=&\frac{\lambda e^{-\lambda t}}{e^{-\lambda x_1}\left(1-e^{-\lambda(x_2-x_1)}\right)}=
\frac{\lambda e^{-\lambda(t-x_1)}}{1-e^{-\lambda(x_2-x_1)}}=f_X(t-x_1,\lambda,x_2-x_1).
\end{eqnarray*}
As we see, the random variable $Y$ has a shifted, truncated distribution with support on an interval of length $x_2-x_1$.

Let the random variable $v_{x,\lambda}$ denote the number of intervals at saturation for this generalized parking process and $M_\lambda(x)=E[v_{x,\lambda}]$. In Section \ref{sec:2} of this paper, we show that $M_{\lambda}(x)$ satisfies 
\begin{align}M_{\lambda}(x+1)&=\displaystyle\int_0^x\frac{\lambda e^{-\lambda t}}{1-e^{-\lambda x}}\left(M_{\lambda}(t)+M_{\lambda}(x-t)\right) \,dt+1\label{eq:3}\\
M_{\lambda}(x)&=0 \mbox{ for } 0\le x\le 1.\nonumber\end{align}
Applying the method of steps to this equation, we can also determine that
\begin{alignat}{3}
M_{\lambda}(x)=1&& &&\quad\mbox{for }1<x\le2,\nonumber\\
M_{\lambda}(x)=1&&\,+\,\frac{(1+e^{-\lambda})(1-e^{-\lambda(x-2)})}{1-e^{-\lambda(x-1)}}&&\quad\mbox{for }2<x\le3\label{eq:18}.
\end{alignat}
Then, letting $\varphi_{\lambda}(s)=\int_0^{\infty}M_{\lambda}(x)e^{-xs}\,dx$ be the Laplace transform of $M_{\lambda}(x)$ and taking the Laplace transform of \eqref{eq:3}, we obtain a functional equation in $\varphi_\lambda(s)$. While we cannot solve this equation explicitly for $\varphi_\lambda(s)$, we can still use Tauberian theory to show that \begin{equation}C_{\lambda}:=\lim_{x\to\infty}\frac{M_\lambda(x)}{x}=\frac{\lambda}{\lambda+1}\left(1+\lambda\varphi_\lambda(\lambda)\right).\label{eq:23}\end{equation} 
This, in turn, can be estimated numerically to arbitrary precision for fixed $\lambda$ by applying the method of steps to equation $\eqref{eq:3}$ and using the estimate \begin{equation}\left\lceil\frac{x-1}{2}\right\rceil\le M_\lambda(x)\le\lfloor x\rfloor,\label{eq:4}\end{equation} which holds for all $x,\lambda>0$. The left inequality is due to the fact that $2\cdot v_{x,\lambda}+1\ge x$ since the gaps between unit intervals at saturation cannot exceed $1$. 

In Section \ref{sec:3}, we use the final value theorem for Laplace transforms to show that
\begin{equation}
B_{\lambda}:=\lim_{x\to\infty}\left(M_\lambda(x)-xC_{\lambda}\right)\label{eq:5}
\end{equation}
exists. We then describe a way to estimate it using \eqref{eq:4} and the method of steps. Unfortunately, as $\lambda\to0^+$, the number of steps needed to achieve a given level of precision for $C_{\lambda}$ and $B_{\lambda}$ increases without bound. We end the section by showing that for each $n=1,2,\ldots$,
\begin{align}
M_\lambda(x) = xC_{\lambda}+B_{\lambda}+o\left(e^{-n x}\right).\label{eq:6}
\end{align}
This result is consistent with the type of decay for the error term described in \eqref{Dvoretsky1}, which could perhaps be proved by generalizing Dvoretsky's methods.

In section \ref{sec:4}, we address the unboundedness problem described above. First, using \eqref{eq:3}, we derive the following integral equation for the derivative of $M_\lambda$:
\begin{align}
M'_{\lambda}(x+1)&=\frac{\int_1^x\lambda\sinh(\lambda t)M'_{\lambda}(t)\,dt+\lambda\sinh(\lambda)}{\cosh(\lambda x)-1}.\label{eq:7}\\
M'_{\lambda}(x)&=0 \mbox{ for } 0< x < 1 \mbox{ and }1<x<2.\nonumber
\end{align}
Then, we generalize an approach of Dvoretsky and Robbins \cite{Dvoretsky} to show that for all $x\ge 3$, 
\begin{equation}
\inf_{x-1\le t\le x}M'_\lambda(t)\le\inf_{x\le t\le x+1}M'_\lambda(t)\le \sup_{x\le t\le x+1}M'_\lambda(t)\le \sup_{x-1\le t\le x}M'_\lambda(t).\label{eq:9}
\end{equation}
Thus, applying numerical methods to \eqref{eq:7} to compute $M_\lambda'(x)$ for small $x$ and $\lambda$ in combination with \eqref{eq:9}, we are able to vastly improve the estimate in \eqref{eq:4}, which in turn allows us to use \eqref{eq:23} to compute $C_\lambda$ and $B_\lambda$ very precisely for small $\lambda>0$.

We should expect that as $\lambda\to0^+$, the results for the exponential distribution converge to R\'{e}nyi 's findings for the case of the uniform distribution. In section \ref{sec:5}, we show that, in fact, $M'_{\lambda}(x)\to M'(x)$ uniformly in $x\in[0,\infty)$ as $\lambda\to0^+$, which implies that 
$$\lim_{\lambda\to0^+}C_\lambda=C\mbox{ and }\lim_{\lambda\to0^+}B_\lambda=B,$$
where $C\approx 0.748$ is the above-mentioned R\'{e}nyi  constant and $B=-(1-C)$ as per \eqref{eq:2}.

In Section \ref{sec:6}, we investigate the asymptotic behavior of $\sigma_{\lambda}^2(x)=\,\, $Var$(v_{x,\lambda})$ by deriving an integral equation for $M_{2,\lambda}=E\left[v_{x,\lambda}^2\right]$. Then, we use Laplace transforms in combination with \eqref{eq:6} to show that the asymptotic slope of the variance,
$$D_{\lambda}=\lim_{x\to\infty}\frac{\sigma_{\lambda}^2(x)}{x}$$
exists and describe its behavior as $\lambda\to\infty$. 
Similar to $M_\lambda(x)$, there exists 
constant $E_\lambda$ such that for each $n=1,2,\ldots$,
\begin{equation}
\sigma_{\lambda}^2(x)=xD_\lambda+E_\lambda+o\left(e^{-n x}\right).\label{40}
\end{equation}
Again, this is consistent with a decay rate similar to \eqref{Dvoretsky2}. We finish the paper by applying the methods from Section \ref{sec:2} to obtain estimates for $D_\lambda$.

In the last section we provide the concluding remarks.
\section{Estimate of asymptotic slope $C_\lambda$}
\label{sec:2}
We begin by showing that \eqref{eq:3} holds. Indeed, let the random variable $T$ denote the left endpoint of the first unit interval placed on the interval $(0,x+1)$, which will divide the original interval $(0,x+1)$ into two disjoint subintervals: $I_1=(0,T)$ and $I_2=(T+1,x+1)$. Let the random variables $v_{x,\lambda}(I_1)$ and $v_{x,\lambda}(I_2)$ denote the number of unit intervals covering $I_1$ and $I_2$, respectively, at saturation. Then, due to the self-similarity of the truncated exponential distribution (shown above), the random variables $v_{x+1,\lambda}(I_1)$ and $v_{T,\lambda}$ have the same distribution, i.e. $v_{x+1,\lambda}(I_1)\sim v_{T,\lambda}$. Similarly, $v_{x+1,\lambda}(I_2)\sim v_{x-T,\lambda}$. Thus, $v_{x+1,\lambda}=v_{T,\lambda}+v_{x-T,\lambda}+1$. Let $M_\lambda(x)=E[v_{x,\lambda}]$. Since the density of the distribution of $T$ is given by $f(t;\lambda,x+1)$, equation \eqref{eq:3} follows.

Now, let 
$$\,\,\,\,\quad\quad\quad\quad\quad\quad\quad\quad\quad\quad\quad\quad\varphi_{\lambda}(s)=\int_0^{\infty}M_{\lambda}(x)e^{-xs}\,dx\quad\quad\quad\quad\quad\quad\quad(s=\sigma+i\omega,\,\,\sigma>0)$$ 
be the Laplace transform of $M_{\lambda}(x)$, which is defined since $M_{\lambda}(x)\le x$ for all $x$. Multiplying \eqref{eq:3} by $1-e^{-\lambda x}$, we obtain
\begin{equation}
M_{\lambda}(x+1)-e^{-\lambda x}M_\lambda(x+1)=1-e^{-\lambda x}+\int_0^x\lambda e^{-\lambda t}M_\lambda(t)\,dt+\int_0^x\lambda e^{-\lambda t}M_\lambda(x-t)\,dt. \label{eq:10}
\end{equation}
Since $M_\lambda(x)=0$ for $0\le x\le 1$, we have 
\begin{align*}
\int_0^\infty M_\lambda(x+1)e^{-xs}\,dx&=\int_1^\infty M_\lambda(u)e^{-(u-1)s}\,du=e^s\int_0^\infty M_\lambda(u)e^{-us}\,du=e^s\varphi_\lambda(s).
\end{align*}
Therefore, taking the Laplace transform of both sides of \eqref{eq:10} yields
\begin{equation}e^s\varphi_\lambda(s)-e^{s+\lambda}\varphi_\lambda(s+\lambda)=\frac{1}{s}-\frac{1}{s+\lambda}+\frac{\lambda\varphi_\lambda(s+\lambda)}{s}+\frac{\lambda\varphi_\lambda(s)}{s+\lambda}.\label{eq:11}
\end{equation}
Multiplying both sides of \eqref{eq:11} by $s^2(s+\lambda)$ and rearranging, we have 
\begin{equation}s^2\varphi_\lambda(s)=\frac{s}{e^s(s+\lambda)-\lambda}\cdot\left(\lambda+(s+\lambda)(se^{s+\lambda}+\lambda)\varphi_\lambda(s+\lambda)\right).\label{eq:12}
\end{equation}
By a Hardy-Littlewood Tauberian theorem\footnote{Let $S:[0,\infty)\to\mathbb R$ be a nondecreasing function such that $S(0)=0$  and let $w(s)=\displaystyle\int_0^\infty e^{-xs}\,dS(x)$. Then, for $\rho>0$, $$w(s)\sim\frac{C}{s^\rho}\quad\mbox{as $s\to0^+$}$$ if and only if $$S(x)\sim\frac{C}{\Gamma(\rho+1)}x^\rho\quad\mbox{as $x\to\infty$}.$$} 
applied to $S(x)=\displaystyle\int_0^xM_\lambda(t)\,dt$ with $\rho=2$, it follows that 
\begin{equation}
\lim_{s\to0^+}s^2\varphi_\lambda(s)=\lim_{x\to\infty}\frac{2\displaystyle\int_0^xM_\lambda(t)\,dt}{x^2}=\lim_{x\to\infty}\frac{M_\lambda(x)}{x}.\label{eq:13}
\end{equation} 
Thus, from \eqref{eq:12} and \eqref{eq:13}, we conclude that \eqref{eq:23} holds, and the problem of estimating $C_\lambda$ is reduced to estimating $\varphi_\lambda(\lambda)$. It follows from \eqref{eq:4} that
\begin{align*}
\int_0^\infty\lambda \left\lfloor x\right\rfloor e^{-\lambda x}\,dx&=\sum_{k=0}^\infty\int_{k}^{k+1}k\lambda e^{-\lambda x}\,dx=\sum_{k=0}^{\infty}\left(ke^{-\lambda k}-ke^{-\lambda(k+1)}\right)\\
&=\sum_{k=0}^{\infty}ke^{-\lambda k}-\sum_{k=1}^{\infty}(k-1)e^{-\lambda k}=\sum_{k=1}^\infty e^{-\lambda k}=\frac{e^{-\lambda}}{1-e^{-\lambda}}\end{align*}
and 
\begin{align*}
\int_0^\infty\lambda \left\lceil \frac{x-1}{2}\right\rceil e^{-\lambda x}\,dx&=\sum_{k=1}^\infty\int_{2k-1}^{2k+1}k\lambda e^{-\lambda x}\,dx=\sum_{k=1}^{\infty}\left(ke^{-\lambda(2k-1)}-ke^{-\lambda(2k+1)}\right)\\
&=\sum_{k=1}^{\infty}ke^{-\lambda(2k-1)}-\sum_{k=2}^{\infty}(k-1)e^{-\lambda(2k-1)}\\
&=e^{-\lambda}+\sum_{k=2}^\infty e^{-\lambda (2k-1)}=e^{-\lambda}+\frac{e^{-3\lambda}}{1-e^{-2\lambda}}=\frac{e^{-\lambda}}{1-e^{-2\lambda}}.
\end{align*}
Therefore, 
\begin{equation}C_\lambda^-:=\frac{\lambda}{\lambda+1}\left(1+\frac{e^{-\lambda}}{1-e^{-2\lambda}}\right)\le C_\lambda\le\frac{\lambda}{\lambda+1}\left(1+\frac{e^{-\lambda}}{1-e^{-\lambda}}\right):=C_\lambda^+.\label{eq:15}\end{equation}
From \eqref{eq:15}, we see that as $\lambda\to\infty$, $$C_\lambda^+-C_\lambda^-=\frac{\lambda e^{-2\lambda}}{(\lambda+1)(1-e^{-2\lambda})}=O\left(e^{-2\lambda}\right)$$ and 
\begin{equation}C_\lambda=\frac{\lambda(1+e^{-\lambda})}{\lambda+1}\,+\,O(e^{-2\lambda}).\label{eq:16}\end{equation}
However, 
\begin{equation}
\lim_{\lambda\to0^+}\left(C_\lambda^+-C_\lambda^-\right)=0.5,\label{eq:17}
\end{equation}
which limits the usefulness of \eqref{eq:15} to large values of $\lambda$. 

Next, we show how to improve \eqref{eq:15}. By using the method of steps applied to equation \eqref{eq:3}, we can numerically compute $M_{\lambda}(x)$ for $0\le x\le n$ for some positive integer $n$. Then, using \eqref{eq:4} for $x>n$, we obtain
\begin{align*}
\int_n^\infty\lambda \left\lfloor x\right\rfloor e^{-x\lambda}\,dx&=\sum_{k=n}^\infty\int_{k}^{k+1}k\lambda e^{-\lambda x}\,dx=\sum_{k=n}^{\infty}\left(ke^{-\lambda k}-ke^{-\lambda(k+1)}\right)\\
&=\sum_{k=n}^{\infty}ke^{-\lambda k}-\sum_{k=n+1}^{\infty}(k-1)e^{-\lambda k}\\
&=ne^{-\lambda n}+\sum_{k=n+1}^\infty e^{-\lambda k}=ne^{-\lambda n}+\frac{e^{-\lambda(n+1)}}{1-e^{-\lambda}}\\
&=\frac{e^{-\lambda n}\left(n-(n-1)e^{-\lambda}\right)}{1-e^{-\lambda}}.\end{align*}
For odd $n$,
\begin{align*}
\int_0^\infty\lambda \left\lceil \frac{x-1}{2}\right\rceil e^{-x\lambda}\,dx&=\sum_{k=n}^\infty\int_{k}^{k+2}\frac{k+1}{2}\lambda e^{-\lambda x}\,dx=\sum_{k=n}^{\infty}\left(\frac{k+1}{2}e^{-\lambda k}-\frac{k+1}{2}e^{-\lambda(k+2)}\right)\\
&=\sum_{k=n}^{\infty}\frac{k+1}{2}e^{-\lambda k}-\sum_{k=n+2}^\infty\frac{k-1}{2}e^{-\lambda k}\\
&=\frac{n+1}{2}e^{-\lambda n}+\frac{n+2}{2}e^{-\lambda(n+1)}+\sum_{k=n+2}^\infty e^{-\lambda k}\\&=\frac{n+1}{2}e^{-\lambda n}+\frac{n+2}{2}e^{-\lambda(n+1)}+\frac{e^{-\lambda(n+2)}}{1-e^{-\lambda}}\\
&=\frac{e^{-\lambda n}\left((n+1)+e^{-\lambda}-(n+1)e^{-2\lambda}\right)}{2(1-e^{-\lambda})},
\end{align*}
and for even $n$,
\begin{align*}
\int_0^\infty\lambda \left\lceil \frac{x-1}{2}\right\rceil e^{-x\lambda}\,dx&=\int_{n}^{n+1}\frac{n}{2}\lambda e^{-\lambda x}\,dx+\sum_{k=n+1}^\infty\int_{k}^{k+1}\frac{k+1}{2}\lambda e^{-\lambda x}\,dx\\
&=\frac{n}{2}e^{-\lambda n}-\frac{n}{2}e^{-\lambda(n+1)}+\sum_{k=n+1}^{\infty}\left(\frac{k+1}{2}e^{-\lambda k}-\frac{k+1}{2}e^{-\lambda(k+2)}\right)\\
&=\frac{n}{2}e^{-\lambda n}-\frac{n}{2}e^{-\lambda(n+1)}+\sum_{k=n+1}^{\infty}\frac{k+1}{2}e^{-\lambda k}-\sum_{k=n+2}^{\infty}\frac{k-1}{2}e^{-\lambda k}\\
&=\frac{n}{2}e^{-\lambda n}-\frac{n}{2}e^{-\lambda(n+1)}+\frac{n+2}{2}e^{-\lambda(n+1)}+\sum_{k=n+2}^{\infty}e^{-\lambda k}\\
&=\frac{n}{2}e^{-\lambda n}+e^{-\lambda(n+1)}+\frac{e^{-\lambda(n+2)}}{1-e^{-\lambda}}=\frac{e^{-\lambda n}\left(n-(n-2)e^{-\lambda}\right)}{2(1-e^{-\lambda})},
\end{align*}
thus giving us the estimate
\begin{equation}C\,^{n^-}_\lambda\le C_\lambda\le C\,^{n^+}_\lambda,\label{101}\end{equation}
where 
$$C\,^{n^-}_\lambda=\begin{cases}\displaystyle\frac{\lambda}{\lambda+1}\left(1+\displaystyle\int_0^n\lambda M_{\lambda}(x)e^{-\lambda x}\,dx+\frac{e^{-\lambda n}}{2(1-e^{-\lambda})}\left((n+1)+e^{-\lambda}-(n+1)e^{-2\lambda}\right)\right)&\quad \mbox{for odd }n\vspace{5mm}\\
\displaystyle\frac{\lambda}{\lambda+1}\left(1+\displaystyle\int_0^n\lambda M_{\lambda}(x)e^{-\lambda x}\,dx+\frac{e^{-\lambda n}}{2(1-e^{-\lambda})}\left(n-(n-2)e^{-\lambda}\right)\right)&\quad \mbox{for even }n\end{cases}$$
and $$C\,^{n^+}_\lambda=\frac{\lambda}{\lambda+1}\left(1+\displaystyle\int_0^n\lambda M_{\lambda}(x)e^{-\lambda x}\,dx+\frac{e^{-\lambda n}}{1-e^{-\lambda}}\left(n-(n-1)e^{-\lambda}\right)\right).$$
Then, 
$$C\,^{n^+}_\lambda-C\,^{n^-}_\lambda=\begin{cases}\displaystyle\frac{\lambda e^{-\lambda n}}{2(\lambda+1)(1-e^{-\lambda})}\left((n-1)e^{-\lambda n}-(2n-1)e^{-\lambda(n+1)}+(n+1)e^{-\lambda(n+2)}\right)&\quad \mbox{for odd }n\vspace{5mm}\\
\displaystyle\frac{\lambda e^{-\lambda n}}{2(\lambda+1)(1-e^{-\lambda})}\left(ne^{-\lambda n}-ne^{-\lambda(n+1)}\right)&\quad \mbox{for even }n.
\end{cases}$$
In either case, $C\,^{n^+}_\lambda-C\,^{n^-}_\lambda=O\left(e^{-2\lambda n}\right)$, but for fixed $n$, 
$$\lim_{\lambda\to0^+}\left(C\,^{n^+}_\lambda-C\,^{n^-}_\lambda\right)=0.5,$$
which is the same as in \eqref{eq:17}. Therefore, this method does indeed lead to an improvement of \eqref{eq:15}, but cannot approximate $C_\lambda$ uniformly for all $\lambda>0$ to an arbitrary level of precision with a fixed number of steps. 
%
%
%
\section{Estimate of the second asymptotic approximation constant $B_\lambda$}
\label{sec:3}
In this section, we improve the estimate for $M_\lambda$ and prove \eqref{eq:6}. Let $f_\lambda(x)=M_\lambda(x)-xC_\lambda$ and $B_\lambda=\displaystyle\lim_{x\to\infty}f_{\lambda}(x)$. We show that $B_\lambda$ exists and describe a way to evaluate it. Let 
$$\quad\quad\quad\quad\quad\quad\quad\quad\hat\varphi_{\lambda}(s)=\int_0^\infty f_\lambda(x)e^{-xs}\,dx=\varphi_\lambda(s)-\frac{C_\lambda}{s^2}\quad\quad\quad\quad\quad\quad(s=\sigma+i\omega,\,\,\sigma>0)$$ 
be the corresponding Laplace transform of $f_{\lambda}(x)$. We seek to apply the final value theorem for Laplace transforms.\footnote{Let $f:[0,\infty)\to\mathbb R$ and let $F(s)=\int_0^\infty f(x)e^{-xs}\,dx$ be its Laplace transform. If $F(s)$ has no poles, $z$, in the complex plane with Re$(z)\ge0$, except at most one pole at $z=0$, then $$\lim_{x\to\infty}f(x)=\lim_{s\to0^+}sF(s).$$}  Indeed, letting $\Phi_\lambda(s)=s^2\varphi_\lambda(s)$, equation \eqref{eq:12} can be written as 
\begin{equation}
\Phi_\lambda(s)=\frac{s}{e^s(s+\lambda)-\lambda}\cdot\left(\lambda+\frac{se^{s+\lambda}+\lambda}{s+\lambda}\cdot\Phi_\lambda(s+\lambda)\right).\label{eq:711}
\end{equation}
Since $\Phi_\lambda(s)$ is clearly analytic for Re$(s)>0$, we see that \eqref{eq:711} defines an analytic continuation of $\Phi_\lambda(s)$ to the half-plane Re$(s)>-\lambda$, and, inductively, to the entire complex plane. Then, in light of \eqref{eq:13}, 
\begin{equation}\varphi_\lambda(s)=\frac{C_\lambda}{s^2}+\frac{B_\lambda}{s}+\sum_{k=0}^\infty a_ks^k\label{eq:19}\end{equation} for all complex $s\ne 0$ for some constants $a_k$. Therefore, $\hat\varphi_\lambda(s)$ is also analytic on the entire complex plane, except for a possible pole of order one at $s=0$, which means the final value theorem applies, and 
\begin{align*}
B_\lambda=\lim_{s\to0^+}s\hat\varphi_\lambda(s)=\lim_{s\to0^+}s\left(\varphi_{\lambda}(s)-\frac{C_\lambda}{s^2}\right)=\lim_{s\to0^+}\frac{s^2\varphi_{\lambda}(s)-C_\lambda}{s}=\lim_{s\to0^+}\left(s^2\varphi_\lambda(s)\right)'.
\end{align*}
From equation \eqref{eq:12}, we see that 
\begin{align*}\left(s^2\varphi_\lambda(s)\right)'=a(\lambda,s)\cdot b(\lambda,s)+ c(\lambda,s)\cdot d(\lambda,s),\end{align*}
where
\begin{align*}a(\lambda,s)&=\frac{e^s(\lambda-\lambda s-s^2)-\lambda}{(e^s(s+\lambda)-\lambda)^2},\\
b(\lambda,s)&=\lambda+(s+\lambda)(se^{s+\lambda}+\lambda)\varphi_\lambda(s+\lambda),\\
c(\lambda,s)&=\frac{s}{e^s(s+\lambda)-\lambda},\\
d(\lambda,s)&=\left(\lambda+e^{s+\lambda}(s^2+(\lambda+2)s+\lambda)\right)\varphi_\lambda(s+\lambda)+(s+\lambda)(se^{s+\lambda}+\lambda)\varphi_\lambda'(s+\lambda).\end{align*}
We have
\begin{align*}
\lim_{s\to0^+}a(\lambda,s)&=\lim_{s\to0^+}\frac{-s(s+\lambda+2)}{2(e^s(s+\lambda)-\lambda)(s+\lambda+1)}=-\frac{\lambda+2}{2(\lambda+1)^2},\\
\lim_{s\to0^+}b(\lambda,s)&=\lambda(1+\lambda\varphi_\lambda(\lambda)),\\
\lim_{s\to0^+}c(\lambda,s)&=\frac{1}{\lambda+1},\\
\lim_{s\to0^+}d(\lambda,s)&=\lambda(1+e^\lambda)\varphi_\lambda(\lambda)+\lambda^2\varphi_\lambda'(\lambda).
\end{align*}
Therefore, \begin{align}B_\lambda&=-\frac{\lambda(\lambda+2)(1+\lambda\varphi_\lambda(\lambda))}{2(\lambda+1)^2}+\frac{\lambda}{\lambda+1}\left((1+e^\lambda)\varphi_\lambda(\lambda)+\lambda\varphi_\lambda'(\lambda)\right)\label{eq:14}\\
&=C_\lambda\left(\frac{2+2e^\lambda+2\lambda e^{\lambda}-\lambda^2}{2\lambda(\lambda+1)}\right)-\frac{e^\lambda+1}{\lambda+1}+\frac{\lambda^2\varphi_\lambda'(\lambda)}{\lambda+1},\nonumber
\end{align}
where the last equality is obtained by substituting the identity in equation \eqref{eq:23} and rearranging. Now, $$\varphi_\lambda'(s)=\frac{d}{ds}\int_0^\infty M_\lambda(x)e^{-xs}\,dx=-\int_0^\infty xM_\lambda(x)e^{-xs}\,dx.$$ Thus, to obtain bounds for $B_\lambda$, we can use the estimate $$x\left\lceil\frac{x-1}{2}\right\rceil\le xM_\lambda(x)\le x\left\lfloor x\right\rfloor,$$
which follows from \eqref{eq:4}, to estimate $\displaystyle\int_0^\infty \lambda^2xM_\lambda(x)e^{-x\lambda}\,dx$. We get
\begin{align*}
\int_0^\infty\lambda^2 x\left\lfloor x\right\rfloor e^{-\lambda x}\,dx&=\sum_{k=0}^\infty\int_{k}^{k+1}k\lambda^2 xe^{-\lambda x}\,dx=\sum_{k=0}^{\infty}\left(k(\lambda k+1)e^{-\lambda k}-k(\lambda(k+1)+1)e^{-\lambda(k+1)}\right)\\
&=\sum_{k=0}^{\infty}k(\lambda k+1)e^{-\lambda k}-\sum_{k=1}^\infty(k-1)(\lambda k+1)e^{-\lambda k}=\sum_{k=1}^{\infty}(\lambda k+1)e^{-\lambda k}\\
&=\lambda\sum_{k=1}^\infty ke^{-\lambda k}+\sum_{k=1}^\infty e^{-\lambda k}=\frac{\lambda}{(1-e^{-\lambda})^2}+\frac{e^{-\lambda}}{1-e^{-\lambda}}\\
&=\frac{e^{-\lambda}(\lambda+1-\lambda e^{-\lambda})}{(1-e^{-\lambda})^2}
\end{align*}
and 
\begin{align*}
\int_0^\infty\lambda^2x\left\lceil \frac{x-1}{2}\right\rceil e^{-\lambda x}\,dx&=\sum_{k=1}^\infty\int_{2k-1}^{2k+1}k\lambda^2 xe^{-\lambda x}\,dx\\
&=\sum_{k=1}^{\infty}\left(k(2\lambda k+1-\lambda)e^{-\lambda(2k-1)}-k(2\lambda k+1+\lambda)e^{-\lambda(2k+1)}\right)\\
&=\sum_{k=1}^{\infty}k(2\lambda k+1-\lambda)e^{-\lambda(2k-1)}-\sum_{k=2}^\infty (k-1)(2\lambda k+1-\lambda)e^{-\lambda(2k-1)}\\
&=(\lambda+1)e^{-\lambda}+\sum_{k=2}^\infty(2\lambda k+1-\lambda)e^{-\lambda(2k-1)}\\
&=(\lambda+1)e^{-\lambda}+2\lambda\sum_{k=2}^\infty ke^{-\lambda(2k-1)}+(1-\lambda)\sum_{k=2}^\infty e^{-\lambda(2k-1)}\\
&=(\lambda+1)e^{-\lambda}+2\lambda\left(\frac{2e^{-3\lambda}}{1-e^{-2\lambda}}+\frac{e^{-5\lambda}}{(1-e^{-2\lambda})^2}\right)+(1-\lambda)\cdot\frac{e^{-3\lambda}}{1-e^{-2\lambda}}\\
&=\frac{e^{-\lambda}(\lambda+1+(\lambda-1)e^{-2\lambda})}{(1-e^{-2\lambda})^2}.
\end{align*}
Thus, 
\begin{equation}B_\lambda^-\le B_\lambda\le B_\lambda^+,\label{eq:26}\end{equation}
where
\begin{align*}B_\lambda^-&=C_\lambda\left(\frac{2+2e^\lambda+2\lambda e^{\lambda}-\lambda^2}{2\lambda(\lambda+1)}\right)-\frac{e^\lambda+1}{\lambda+1}-\frac{e^{-\lambda}(\lambda+1-\lambda e^{-\lambda})}{(\lambda+1)(1-e^{-\lambda})^2},\\
B_\lambda^+&=C_\lambda\left(\frac{2+2e^\lambda+2\lambda e^{\lambda}-\lambda^2}{2\lambda(\lambda+1)}\right)-\frac{e^\lambda+1}{\lambda+1}-\frac{e^{-\lambda}(\lambda+1+(\lambda-1)e^{-2\lambda})}{(\lambda+1)(1-e^{-2\lambda})^2}.\end{align*}
As before, this estimate is only useful for large values of $\lambda$. Finally, it follows from \eqref{eq:15} that
\begin{equation}B_\lambda=-\frac{1}{2}+\frac{1}{\lambda+1}+\frac{1}{2(\lambda+1)^2}+O\left(e^{-\lambda}\right).\label{eq:21}\end{equation}
Next, we verify \eqref{eq:6}, which is a special case of the following result.

\textbf{Lemma 1:} Let $f:[0,\infty)\mapsto\mathbb R$. Suppose that its Laplace transform $F(s)=\int_0^\infty f(x)e^{-xs}\,dx$ is analytic on the entire complex plane. Then, for each $n=1,2,\ldots$, $f(x)=o(e^{-nx})$.

\noindent\underline{Proof}: For each $n=1,2,\ldots$, let $f_n(x)=e^{nx}f(x)$. Then, its Laplace transform $F_n(s)=F(s-n)$ is also analytic on the entire complex plane, so $\displaystyle\lim_{s\to0^+}sF_n(s)=0$, and Lemma 1 follows by the final value theorem.\qed
 \mbox{  }\\
Thus, \eqref{eq:6} follows from  \eqref{eq:19} and Lemma 1 applied to  the function $M_\lambda(x)-xC_\lambda-B_\lambda$.\\

\textbf{Remark 1:} R\'{e}nyi  \cite{Renyi} gives an intuitive justification for \eqref{eq:2} based on the fact that the function $f(x)=Cx-(1-C)$ satisfies the integral equation in \eqref{eq:1} for any constant $C$. It is noteworthy that $f(x)$ also satisfies the equation in \eqref{eq:3} for any constants $C$ and $\lambda$, but the analogous result does not hold. In other words, $B_\lambda\ne -(1-C_\lambda)$ for arbitrary $\lambda$, which can be seen from \eqref{eq:21} and from the figures in the next section.
\section{Improvement of estimates for $C_\lambda$ and $B_\lambda$.}
\label{sec:4}
In order to improve the accuracy of the estimates described in the previous sections, we now take a slightly different approach. First, we show that \eqref{eq:7} holds. Indeed, differentiating both sides of \eqref{eq:10} with respect to $x$, we obtain
$$\lambda e^{-\lambda x}M_{\lambda}(x+1)+(1-e^{-\lambda x})M'_{\lambda}(x+1)=\lambda e^{-\lambda x}+\lambda e^{-\lambda x}M_{\lambda}(x)+\lambda M_{\lambda}(x)-\lambda^2\int_0^x e^{-\lambda (x-t)}M_{\lambda}(t)\,dt.$$
Multiplying by $e^{\lambda x}$, differentiating once more, and simplifying, we get
$$\lambda M'_{\lambda}(x+1)+\lambda e^{\lambda x}M'_{\lambda}(x+1)+(e^{\lambda x}-1)M''_{\lambda}(x+1)=\lambda (e^{\lambda x}+1)M'_{\lambda}(x).$$
Then, multiplying by the integrating factor $\frac{e^{-\lambda x/2}}{2}\sinh\left(\frac{\lambda x}{2}\right)$ gives us 
$$\left(\sinh^2\left(\frac{\lambda x}{2}\right)M'_{\lambda}(x+1)\right)'=
\frac{\lambda}{2}\sinh(\lambda x)M'_{\lambda}(x).$$
Since $M'_\lambda(x)$ is continuous for $x\ge 1$, integrating and simplifying yields
\begin{equation}
\left(\cosh(\lambda x)-1\right)M'_{\lambda}(x+1)=\int_1^x\lambda\sinh(\lambda t)M'_{\lambda}(t)\,dt+\left(\cosh(\lambda)-1\right)M'_{\lambda}(2^+)\label{eq:22}.
\end{equation} 
It remains to determine $M'_{\lambda}(2^+)$. Differentiating \eqref{eq:18} and simplifying, we see that for \\$2<x\le 3$,
$$M'_\lambda(x)=\frac{\lambda\sinh(\lambda)}{\cosh(\lambda(x-1))-1},$$
so $$M'_\lambda(2^+)=\frac{\lambda\sinh(\lambda)}{\cosh(\lambda)-1}.$$
Substituting this into \eqref{eq:22} and dividing both sides by $\left(\cosh(\lambda x)-1\right)$ yields \eqref{eq:7}.

Next, we show that  \eqref{eq:9} holds for $x\ge3$. This is a consequence of the following theorem, which is similar to Theorem $1$ of Dvoretsky and Robbins \cite{Dvoretsky}.
\textbf{Theorem 1:} Let $K(t)\ge0$ for all $t>0$, and suppose that for all $x\ge a\ge0$, $f(x)$ satisfies 
\begin{equation}f(x+1)=\frac{\int_0^xK(t)f(t)\,dt+C}{\int_0^xK(t)\,dt}.\label{eq:24}\end{equation}
Then, for all $x\ge a+2$, 
\begin{equation}
I_{x+1}\le I_x\le S_x\le S_{x+1},\label{eq:25}
\end{equation}
where $I_x=\displaystyle\inf_{x-1\le t\le x}f(t)$ and $S_x=\displaystyle\sup_{x-1\le t\le x}f(t)$.

\noindent \underline{Proof}: For any $x\ge a+2$ and any $t\in[x,x+1]$, we have
\begin{align*}f(t)&=\frac{\int_0^{t-1}K(s)f(s)\,ds+C}{\int_0^{t-1}K(s)\,ds}=\frac{\int_0^{x-1}K(s)f(s)\,ds+C}{\int_0^{t-1}K(s)\,ds}+\frac{\int_{x-1}^{t-1}K(s)f(s)\,ds}{\int_0^{t-1}K(s)\,ds}\\
&=\frac{\int_0^{x-1}K(s)\,ds}{\int_0^{t-1}K(s)\,ds}\cdot f(x)+\frac{\int_{x-1}^{t-1}K(s)f(s)\,ds}{\int_0^{t-1}K(s)\,ds}\\
&\ge\frac{\int_0^{x-1}K(s)\,ds}{\int_0^{t-1}K(s)\,ds}\cdot I_x+\frac{\int_{x-1}^{t-1}K(s)I_x\,ds}{\int_0^{t-1}K(s)\,ds}=I_x.
\end{align*}
Thus, the first inequality in \eqref{eq:25} holds. The third inequality is established in the same way by replacing ``$\ge$" with ``$\le$" and ``$\inf$" with ``$\sup$" in the above proof.\qed

Since $M'_{\lambda}(x)=0$ for $0<x<1$, \eqref{eq:7} implies that $$M'_{\lambda}(x+1)=\frac{\int_0^x\lambda\sinh(\lambda t)M'_{\lambda}(t)\,dt+\lambda\sinh(\lambda)}{\cosh(\lambda x)-1}$$ holds for $x\ge1$, so by Theorem $1$, \eqref{eq:9} holds for $x\ge3$, as desired.

We can use this result to greatly improve our estimates of $B_\lambda$ and $C_\lambda$ in the following way. First, we apply numerical methods to \eqref{eq:7} to compute $I_n$ and $S_n$ for some integer $n$. Then, we conclude inductively from \eqref{eq:9} that for all $x\ge n$, $I_n\le M'_{\lambda}(x)\le S_n$, which implies that
\begin{equation}M_{\lambda}(n)+I_n(x-n)\le M_{\lambda}(x)\le M_{\lambda}(n)+S_n(x-n).\label{100}\end{equation}
Finally, computing $M_{\lambda}(x)$ numerically for $0\le x\le n$ and substituting into \eqref{eq:23} and \eqref{eq:26} leads to very accurate estimates for $C_\lambda$ and $B_\lambda$, respectively, as shown in the graphs above.\\

\begin{figure}
\includegraphics[scale=0.19]{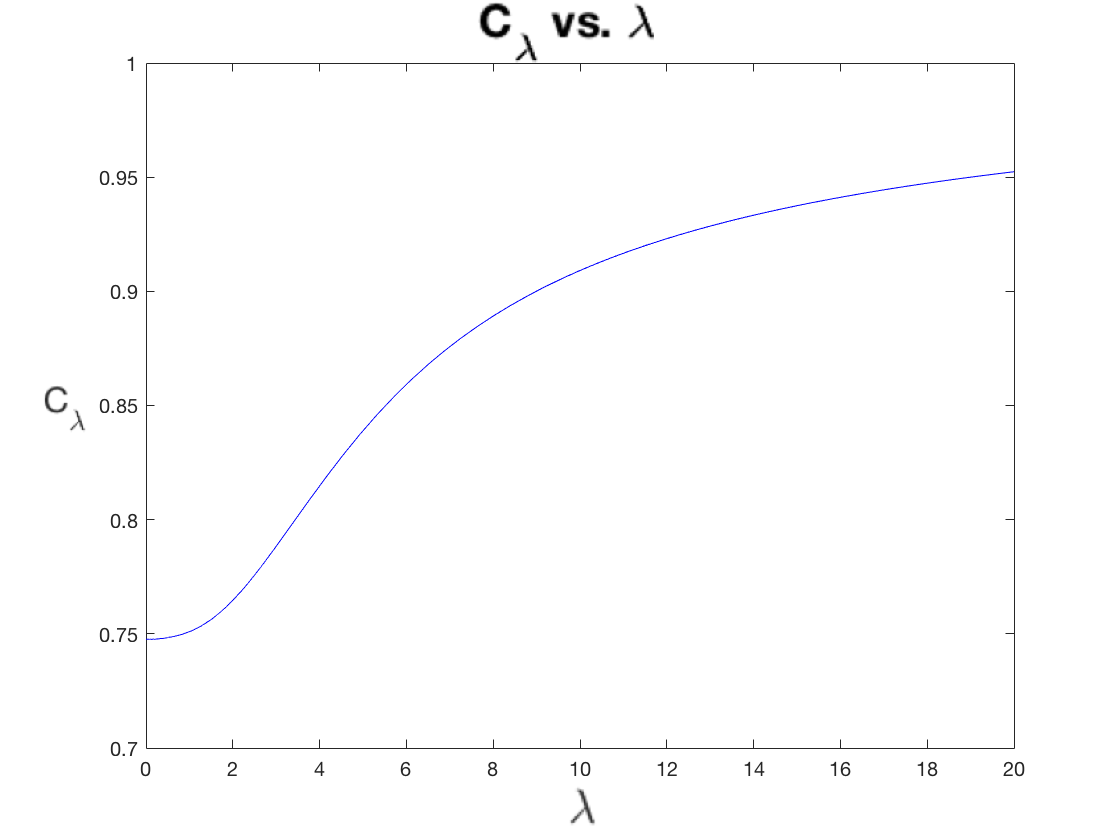}\includegraphics[scale=0.19]{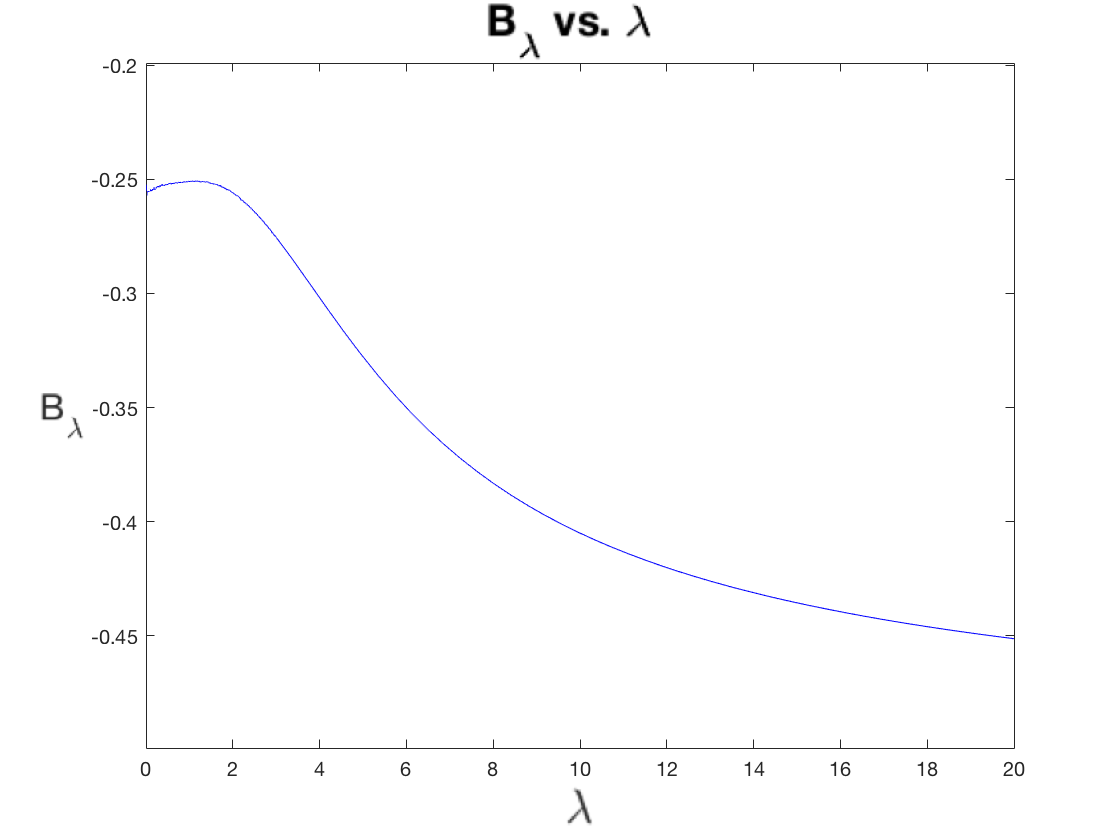}\label{fig:1}\\ 
\caption{The graphs of $C_\lambda$ and $B_\lambda$ vs. $\lambda$. Values for $0<\lambda<3$ were computed using the estimate in \eqref{100} with $n=7$ combined with \eqref{eq:23} and \eqref{eq:26}. For $\lambda\ge3$, computations were done using the estimates in \eqref{101} with $n=7$ and \eqref{eq:26}. The results are consistent with the asymptotic behavior described in \eqref{eq:16} and \eqref{eq:21} as $\lambda\to\infty$, as well as with R\'{e}nyi 's findings as $\lambda\to0^+$, which will be proved in the next section.}
\end{figure}

\textbf{Remark 2:} By an inductive argument applied to \eqref{eq:7}, one can easily show that for fixed $n\ge 3$, as $\lambda\to\infty,$ $$\inf_{n-1\le t\le n}M_{\lambda}'(t)\to0^+\,\,\mbox{  and  }\,\,\sup_{n-1\le t\le n}M_{\lambda}'(t)\to\infty.$$ Thus, the estimate in \eqref{eq:9} suffers from an unboundedness problem similar to the ones described at the ends of Sections \ref{sec:2} and \ref{sec:3} for $\lambda\to0^+$. Thus, we see that the method presented in this section complements the ones from the previous sections, allowing us to estimate $C_\lambda$ and $B_\lambda$ for any values of $\lambda>0$ to arbitrary precision. 

\section{Convergence to uniform case as $\lambda\to 0$.}
\label{sec:5}
In this section, we prove that $M_{\lambda}'(x)\to M'(x)$ uniformly as $\lambda\to0^+$, where $M(x)$ satisfies \eqref{eq:1}. 
It will follow that $\displaystyle\lim_{\lambda\to0^+}C_\lambda=C$ and $\displaystyle\lim_{\lambda\to0^+}B_\lambda=B=-(1-C)$. 

First, we show that $M'(x)$ satisfies
\begin{align}M'(x+1)&=\frac{\int_{1}^x2tM'(t)\,dt+2}{x^2}.\label{eq:27}\\
M'(x)&=0 \mbox{ for } 0< x < 1 \mbox{ and }1<x<2.\nonumber\end{align}
Indeed, multiplying \eqref{eq:1} by $x$ and differentiating twice, we obtain 
$$2M'(x+1)+xM''(x+1)=2M'(x).$$
Multiplying both sides by the integrating factor $x$ gives us
$$\left(x^2M'(x+1)\right)'=2M'(x).$$
\eqref{eq:27} then follows by integrating and substituting initial conditions, similarly to what was done above in the derivation of \eqref{eq:7}.

Next, we prove the following result.

\textbf{Lemma 2:} $M_{\lambda}'(x)\to M'(x)$ uniformly as $\lambda\to0^+$ for all $x$ on compact subsets of $[0,\infty)$.

\noindent\underline{Proof}: We use induction on $n$ to show that the desired convergence holds for all $x\in[0,n]$. This is trivial for $n=1$ and $n=2$, so we assume that $n\ge 2$. Then, for any $x\in[n,n+1]$, 
\begin{align*}
&\left|M_{\lambda}'(x+1)-M'(x+1)\right|=\left|\frac{\int_1^x\lambda\sinh(\lambda t)M'_{\lambda}(t)\,dt+\lambda\sinh(\lambda)}{\cosh(\lambda x)-1}-\frac{\int_{1}^x2tM'(t)\,dt+2}{x^2}\right|\\
&=\left|\int_1^x\left(\frac{\lambda\sinh(\lambda t)}{\cosh(\lambda x)-1}-\frac{2t}{x^2}\right)M_{\lambda}'(t)\,dt+\int_1^x\frac{2t}{x^2}\left(M'(t)-M'_\lambda(t)\right)\,dt+\frac{\lambda\sinh(\lambda)}{\cosh(\lambda x)-1}-\frac{2}{x^2}\right|\\
&\le\int_1^x\left|\frac{\lambda\sinh(\lambda t)}{\cosh(\lambda x)-1}-\frac{2t}{x^2}\right|M_{\lambda}'(t)\,dt+\int_1^x\frac{2t}{x^2}\left|M'(t)-M'_\lambda(t)\right|\,dt+\left|\frac{\lambda\sinh(\lambda)}{\cosh(\lambda x)-1}-\frac{2}{x^2}\right|.
\end{align*}
Now, it is readily verified that $$\lim_{\lambda\to0^+}\left(\frac{\lambda\sinh(\lambda t)}{\cosh(\lambda x)-1}-\frac{2t}{x}\right)=0\quad\mbox{ and }\quad\lim_{\lambda\to0^+}\left(\frac{\lambda\sinh(\lambda)}{\cosh(\lambda x)-1}-\frac{2}{x^2}\right)=0$$
for all $2\le t\le x\le n$. Moreover, the convergence is uniform in each case because the respective derivatives in $x$ (and $t$) are bounded on $[2,n]$ for any fixed $n\ge 2$. Thus, for any $\varepsilon>0$, we may choose $\delta_1>0$ such that for all $\lambda\in(0,\delta_1)$,
$$\left|\frac{\lambda\sinh(\lambda t)}{\cosh(\lambda x)-1}-\frac{2t}{x}\right|<\frac{\varepsilon}{3n\overline{M}}\quad\mbox{ and }\quad\left|\frac{\lambda\sinh(\lambda)}{\cosh(\lambda x)-1}-\frac{2}{x^2}\right|<\frac{\varepsilon}{3},$$
where $2\le t\le x\le n$ and $\overline{M}=\displaystyle\max_{2\le t\le n}M'_{\lambda}(t)$.
Moreover, by the inductive hypothesis, we may choose $\delta_2>0$, such that for all $\lambda\in(0,\delta_2)$, $\left|M'_\lambda(t)-M'(t)\right|<\frac{\varepsilon}{6}$ for all $t\in[0,n]$. Therefore, letting $\delta=\min\{\delta_1,\delta_2\}$, we have that for all $x\in[n-1,n]$, $\lambda\in(0,\delta)$ implies that
\begin{align*}
\left|M_{\lambda}'(x+1)-M'(x+1)\right|&\le\int_1^x\left|\frac{\lambda\sinh(\lambda t)}{\cosh(\lambda x)-1}-\frac{2t}{x^2}\right|M_{\lambda}'(t)\,dt+\int_1^x\frac{2t}{x^2}\left|M'(t)-M'_\lambda(t)\right|\,dt\\
&\quad\quad+\left|\frac{\lambda\sinh(\lambda)}{\cosh(\lambda x)-1}-\frac{2}{x^2}\right|\\
&< x\cdot\frac{\varepsilon}{3n}+\frac{\varepsilon}{3}+\frac{\varepsilon}{3}<\varepsilon.
\end{align*}
Thus, $M'_{\lambda}(x)\to M(x)$ as $\lambda\to0^+$ uniformly on $[n,n+1]$ and therefore on $[0,n+1]$ as well.\qed

\textbf{Theorem 2:} $M'_{\lambda}(x)\to M'(x)$ uniformly as $\lambda\to0^+$ for $x$ on $(0,\infty)$. 

\noindent\underline{Proof}: Let $\varepsilon>0$ be given. Since $\displaystyle\lim_{x\to\infty}M'(x)$ exists, we can find $x'\ge1$ such that $$\sup_{x'-1\le t\le x'}M'(t)-\inf_{x'-1\le t\le x'}M'(t)<\varepsilon.$$
It follows from Lemma 2 that there exists $\delta>0$ such that for all $\lambda\in(0,\delta)$, \\$|M'_\lambda(x)-M'(x)|\,\,<\varepsilon$ for all $t\in[0,x']$. Then,
\begin{align*}
\sup_{x'\le t\le x'+1}M'_\lambda(t)&\le\sup_{x'-1\le t\le x'}M'_\lambda(t)<\inf_{x'-1\le t\le x'}M'(t)+\varepsilon\le\inf_{x'\le t\le x'+1}M'(t)+\varepsilon,\\
\inf_{x'\le t\le x'+1}M'_\lambda(t)&\ge\inf_{x'-1\le t\le x'}M'_\lambda(t)>\sup_{x'-1\le t\le x'}M'(t)-\varepsilon\ge\sup_{x'\le t\le x'+1}M'(t)-\varepsilon.
\end{align*}
Thus, $|M'_\lambda(x)-M'(x)|\,\,<\varepsilon$ for all $x\in[x',x'+1]$. It follows inductively that this inequality holds for all $x\ge0$.\qed\vspace{2mm}
\section{Estimate for the variance}
\label{sec:6}
Let $t$ denote the left endpoint of the first unit interval placed on the interval $(0,x+1)$. By an argument similar to the one presented at the beginning of section \ref{sec:2}, $v_{x+1,\lambda}=1+v_{t,\lambda}+v_{x-t,\lambda}$, so 
\begin{align*}
v_{x+1,\lambda}^2&=\left(1+v_{t,\lambda}+v_{x-t,\lambda}\right)^2=1+v_{t,\lambda}^2+v_{x-t,\lambda}^2+2v_{t,\lambda}+2v_{x-t,\lambda}+2v_{t,\lambda}v_{x-t,\lambda}.
\end{align*}
Thus, $M_{2,\lambda}(x):=E\left[v_{x,\lambda}^2\right]$ satisfies $M_{2,\lambda}(x)=0$ for $0\le x\le 1$ and 
$$M_{2,\lambda}(x+1)=\int_0^x \frac{\lambda e^{-\lambda t}}{1-e^{-\lambda x}}\left(M_{2,\lambda}(t)+M_{2,\lambda}(x-t)+2M_{\lambda}(t)+2M_{\lambda}(x-t)+2M_\lambda(t)M_{\lambda}(x-t)\right)\,dt+1.$$
Multiplying by $1-e^{-\lambda x}$ on both sides and simplifying, we obtain
\begin{align}M_{2,\lambda}(x+1)-e^{-\lambda x}M_{2,\lambda}(x+1)&=1-e^{-\lambda x}+\int_0^xM_{2,\lambda}(t)\lambda e^{-\lambda t}\,dt\label{eq:28}\\
&\quad+\int_0^xM_{2,\lambda}(x-t)\lambda e^{-\lambda t}\,dt+2\int_0^xM_{\lambda}(t)\lambda e^{-\lambda t}\,dt\nonumber\\
&\quad+2\int_0^xM_{\lambda}(x-t)\lambda e^{-\lambda t}\,dt\nonumber\\
&\quad+2\int_0^xM_{\lambda}(t)M_\lambda(x-t)\lambda e^{-\lambda t}\,dt\nonumber.
\end{align}
Let $$\varphi_{2,\lambda}(s)=\int_0^\infty M_{2,\lambda}(x)e^{-xs}\,dx$$
be the Laplace transform of $M_{2,\lambda}(x)$. Taking the Laplace transform of both sides of \eqref{eq:28}, we get
\begin{align}
e^s\varphi_{2,\lambda}(s)-e^{s+\lambda}\varphi_{2,\lambda}(\lambda+s)&=\frac{1}{s}-\frac{1}{\lambda+s}+\frac{\lambda \varphi_{2,\lambda}(s+\lambda)}{s}+\frac{\lambda \varphi_{2,\lambda}(s)}{s+\lambda}+\frac{2\lambda \varphi_{\lambda}(s+\lambda)}{s}\\
&\quad+\frac{2\lambda \varphi_{\lambda}(s)}{s+\lambda}+2\lambda\varphi_{\lambda}(s)\varphi_{\lambda}(s+\lambda)\nonumber.
\end{align}
Multiplying by $s^3(s+\lambda)$ and rearranging yields
\begin{align}
s^3\varphi_{2,\lambda}(s)&=\frac{s}{e^s(s+\lambda)-\lambda}\left(\lambda s+\lambda s(s+\lambda)\varphi_{2,\lambda}(s+\lambda)+2\lambda s(s+\lambda)\varphi_\lambda(s+\lambda)\right.\label{eq:29}\\
&\quad\quad\quad\quad\quad\quad\quad+2\lambda s^2\varphi_\lambda(s)+2\lambda s^2(s+\lambda)\varphi_\lambda(s)\varphi_\lambda(s+\lambda)\nonumber\\
&\quad\quad\quad\quad\quad\quad\quad\left.+\,e^{s+\lambda}(s+\lambda)s^2\varphi_{2,\lambda}(s+\lambda)\right)\nonumber.
\end{align}
Thus, since $\displaystyle\lim_{s\to0^+}s^2\varphi_{\lambda}(s)=C_\lambda$, we have
\begin{align*}
\lim_{s\to0^+}s^3\varphi_{2,\lambda}(s)&=\frac{2\lambda C_\lambda+2\lambda^2C_\lambda\varphi_\lambda(\lambda)}{\lambda+1}=\frac{2\lambda C_\lambda+2C_\lambda\left(\lambda C_\lambda+C_\lambda-\lambda\right)}{\lambda+1}=2C_\lambda^2.
\end{align*}
It follows from the previously mentioned Tauberian theorem$^1$ that 
$$\lim_{x\to\infty}\frac{\varphi_{2,\lambda}(x)}{x^2}=C_\lambda^2.$$
By a similar, albeit slightly more tedious, argument to the one presented in Section \ref{sec:3} for $M_\lambda(x)$, \eqref{eq:29} and Lemma $1$ imply that there exist constants, $B_{1,\lambda}$ and $B_{0,\lambda}$, such that for each $n=1,2,\ldots$, 
\begin{equation}M_{2,\lambda}(x)=x^2C_\lambda^2+x B_{1,\lambda}+B_{0,\lambda}+o\left(e^{-n x}\right)\label{eq:30}.\end{equation}
From \eqref{eq:6} and \eqref{eq:30}, we conclude that \eqref{40} holds with 
\begin{equation}
D_\lambda=B_{1,\lambda}-2B_\lambda C_\lambda.\label{30} 
\end{equation} 
Next, we examine the asymptotic behavior of $D_\lambda$ as $\lambda\to\infty$. In light of \eqref{30}, since this was already done for $B_\lambda$ and $C_\lambda$ in the previous sections, it remains to do the same for $B_{1,\lambda}$. Let $h_{\lambda}(x)=M_{2,\lambda}(x)-x^2C_\lambda^2$ and 
\begin{equation}\tilde{\varphi}_\lambda(s)=\int_0^\infty h_{\lambda}(x)e^{-xs}\,dx=\varphi_{2,\lambda}(s)-\frac{2C_\lambda^2}{s^3}\nonumber\end{equation}
be its Laplace transform. Then, 
\begin{equation}D_\lambda=\lim_{s\to0^+}s^2\tilde{\varphi}_\lambda(s)=\lim_{s\to0^+}\frac{s^3\varphi_{2,\lambda}(s)-2C_\lambda^2}{s}=\lim_{s\to0^+}\left(s^3\varphi_{2,\lambda}(s)\right)'.\end{equation}
From equation \eqref{eq:29}, we see that 
\begin{align*}\left(s^3\varphi_\lambda(s)\right)'=e(\lambda,s)\cdot f(\lambda,s)+ g(\lambda,s)\cdot h(\lambda,s),\end{align*}
where
\begin{align*}e(\lambda,s)&=\frac{e^s(\lambda-\lambda s-s^2)-\lambda}{(e^s(s+\lambda)-\lambda)^2},\\
f(\lambda,s)&=\lambda s+\lambda s(s+\lambda)\varphi_{2,\lambda}(s+\lambda)+2\lambda s(s+\lambda)\varphi_\lambda(s+\lambda)+2\lambda s^2\varphi_\lambda(s)\\
&\quad+2\lambda s^2(s+\lambda)\varphi_\lambda(s)\varphi_\lambda(s+\lambda)+\,e^{s+\lambda}(s+\lambda)s^2\varphi_{2,\lambda}(s+\lambda),\\
g(\lambda,s)&=\frac{s}{e^s(s+\lambda)-\lambda},\\
h(\lambda,s)&=\lambda+\lambda (2s+\lambda)\varphi_{2,\lambda}(s+\lambda)+\lambda s(s+\lambda)\varphi'_{2,\lambda}(s+\lambda)+2\lambda (2s+\lambda)\varphi_\lambda(s+\lambda)\\
&\quad+2\lambda s(s+\lambda)\varphi'_\lambda(s+\lambda)+2\lambda\left(s^2\varphi_\lambda(s)\right)'+2\lambda s^2\varphi_\lambda(s)\varphi_\lambda(s+\lambda)\\
&\quad+2\lambda(s+\lambda)\left(s^2\varphi_\lambda(s)\right)'\varphi_\lambda(s+\lambda)+2\lambda s^2(s+\lambda)\varphi_\lambda(s)\varphi'_\lambda(s+\lambda)\\
&\quad +e^{s+\lambda}\left((s+\lambda)s^2\varphi_{2,\lambda}(s+\lambda)+(3s^2+2\lambda s)\varphi_{2,\lambda}(s+\lambda)+(s+\lambda)s^2\varphi'_{2,\lambda}(s+\lambda)\right).\end{align*}
Using \eqref{eq:23} and \eqref{eq:14} to simplify, we end up with
\begin{align*}
\lim_{s\to0^+}e(\lambda,s)&=-\frac{\lambda+2}{2(\lambda+1)^2},\\
\lim_{s\to0^+}f(\lambda,s)&=2\lambda C_\lambda+2\lambda^2C_\lambda\varphi_\lambda(\lambda)=2(\lambda+1)C_\lambda^2,\\
\lim_{s\to0^+}g(\lambda,s)&=\frac{1}{\lambda+1},\\
\lim_{s\to0^+}h(\lambda,s)&=\lambda+\lambda^2\varphi_{2,\lambda}(\lambda)+2\lambda^2\varphi_\lambda(\lambda)+2\lambda B_\lambda+2\lambda C_\lambda\varphi_\lambda(\lambda)+2\lambda^2B_\lambda\varphi_\lambda(\lambda)+2\lambda^2 C_\lambda\varphi'_\lambda(\lambda)\\
&=\lambda+\lambda^2\varphi_{2,\lambda}(\lambda)+2\left(\lambda C_\lambda+C_\lambda-\lambda\right)+2\lambda B_\lambda+2C_\lambda\left(C_\lambda+\frac{C_\lambda}{\lambda}-1\right)\\
&\quad+2B_\lambda\left(\lambda C_\lambda+C_\lambda-\lambda\right)+2C_\lambda\left(B_\lambda(\lambda+1)-C_\lambda\left(\frac{2+2e^\lambda+2\lambda e^\lambda-\lambda^2}{2\lambda}\right)+e^\lambda+1\right)\\
&=-\lambda+\lambda^2\varphi_{2,\lambda}(\lambda)+2\left(1+\lambda+e^\lambda\right)C_\lambda+4(\lambda+1)B_\lambda C_\lambda+\left(\frac{2\lambda-2e^\lambda-2\lambda e^\lambda+\lambda^2}{\lambda}\right)C_\lambda^2.
\end{align*}
Therefore, 
\begin{equation}B_{1,\lambda}=4 B_\lambda C_\lambda-\left(\frac{2e^\lambda}{\lambda}\right)C_\lambda^2+\left(\frac{2(1+\lambda+e^\lambda)}{\lambda+1}\right)C_\lambda+\frac{\lambda^2\varphi_{2,\lambda}(\lambda)-\lambda}{\lambda+1}.\label{34}\end{equation}
Thus, to obtain bounds for $B_{1,\lambda}$, we can use the estimate 
$$\left\lceil\frac{x-1}{2}\right\rceil^2\le M_{2,\lambda}(x)\le \left\lfloor x\right\rfloor^2$$
to approximate $\lambda\varphi_{2,\lambda}(\lambda)=\displaystyle\int_0^\infty \lambda M_{2,\lambda}(x)e^{-x\lambda}\,dx$. We get
\begin{align*}
\int_0^\infty\lambda \left\lfloor x\right\rfloor^2 e^{-\lambda x}\,dx&=\sum_{k=0}^\infty\int_{k}^{k+1}\lambda k^2e^{-\lambda x}\,dx=\sum_{k=0}^{\infty}\left(k^2e^{-\lambda k}-k^2e^{-\lambda(k+1)}\right)\\
&=\sum_{k=0}^{\infty}k^2e^{-\lambda k}-\sum_{k=1}^\infty(k-1)^2e^{-\lambda k}=\sum_{k=1}^{\infty}(2k-1)e^{-\lambda k}\\
&=\frac{2e^{-\lambda}}{(1-e^{-\lambda})^2}+\frac{e^{-\lambda}}{1-e^{-\lambda}}=\frac{e^{-\lambda}(1+e^{-\lambda})}{(1-e^{-\lambda})^2}
\end{align*}
and 
\begin{align*}
\int_0^\infty\lambda\left\lceil \frac{x-1}{2}\right\rceil^2e^{-\lambda x}\,dx&=\sum_{k=1}^\infty\int_{2k-1}^{2k+1}\lambda k^2e^{-\lambda x}\,dx=\sum_{k=1}^{\infty}\left(k^2e^{-\lambda(2k-1)}-k^2e^{-\lambda(2k+1)}\right)\\
&=\sum_{k=1}^{\infty}k^2e^{-\lambda(2k-1)}-\sum_{k=2}^\infty (k-1)^2e^{-\lambda(2k-1)}=e^{-\lambda} + \sum_{k=2}^\infty(2k-1)e^{-\lambda(2k-1)}\\
&=e^{-\lambda}+\frac{e^{-3\lambda}(3-e^{-2\lambda})}{(1-e^{-2\lambda})^2}=\frac{e^{-\lambda}(1+e^{-2\lambda})}{(1-e^{-2\lambda})^2}.
\end{align*}
Thus, from \eqref{30} and \eqref{34}, we obtain
\begin{equation}D_{\lambda}^-\le D_{\lambda}\le D_{\lambda}^+,\label{eq:35}\end{equation}
where
\begin{align*}D_\lambda^-&=2 B_\lambda C_\lambda-\left(\frac{2e^\lambda}{\lambda}\right)C_\lambda^2+\left(\frac{2(1+\lambda+e^\lambda)}{\lambda+1}\right)C_\lambda+\frac{\lambda}{\lambda+1}\left(\frac{e^{-\lambda}(1+e^{-2\lambda})}{(1-e^{-2\lambda})^2}-1\right),\\
D_\lambda^+&=2 B_\lambda C_\lambda-\left(\frac{2e^\lambda}{\lambda}\right)C_\lambda^2+\left(\frac{2(1+\lambda+e^\lambda)}{\lambda+1}\right)C_\lambda+\frac{\lambda}{\lambda+1}\left(\frac{e^{-\lambda}(1+e^{-\lambda})}{(1-e^{-\lambda})^2}-1\right).\end{align*}
As before, this estimate is only useful for large values of $\lambda$. Finally, from \eqref{eq:16} and \eqref{eq:21}, we conclude that
\begin{equation*}D_\lambda=\frac{\lambda}{(\lambda+1)^3}+O\left(e^{-\lambda}\right) \mbox{ as }\lambda\to\infty.\end{equation*}\\

\textbf{Remark 3:} Dvoretsky and Robbins \cite{Dvoretsky} show that the random variable 
$$Z_x=\frac{v_x-M(x)}{\sigma(x)}\sim N(0,1)$$
asymptotically as $x\to\infty$, where $\sigma(x)$ denotes the standard deviation of $v_x$. They provide two proofs of this result, the second of which is based solely on the relation $\sigma^2(x) = cx+o(x)$, where $c>0$ is constant. By a similar argument to the one presented in their paper, this relation also holds for $\sigma_\lambda(x)$ for each $\lambda>0$. Therefore, we can conclude that the distribution of the random variable 
$$Z_{x,\lambda}=\frac{v_{x,\lambda}-M_\lambda(x)}{\sigma_\lambda(x)}\sim N(0,1)$$
tends to the standard normal distribution asymptotically as $x\to\infty$ for every fixed $\lambda$.

\section{Conclusion} 
\label{conclusion}
By exploiting the self-similarity of the exponential distribution, we analyze the generalized parking problem and derive a functional equation for the mathematical expectation of the corresponding random variable denoting the number of unit intervals that have been placed on the interval $(0,x)$ at saturation. Then, using Laplace transforms and Tauberian theory, we were able to prove the following relations for the expectation and variance for each $n=1,2,\ldots$ and all $\lambda>0$: $$M_\lambda(x) = xC_{\lambda}+B_{\lambda}+o\left(e^{-n x}\right)\,\,\mbox{ and }\,\,\mbox{Var}(v_{x,\lambda})=xD_\lambda+E_\lambda+o\left(e^{-n x}\right),$$
where 
\begin{align*}
C_\lambda&=\frac{\lambda}{\lambda+1}\left(1+e^{-\lambda}\right)+O\left(e^{-2\lambda}\right),\\
B_\lambda&=-\frac{1}{2}+\frac{1}{\lambda+1}+\frac{1}{2(\lambda+1)^2}+O\left(e^{-\lambda}\right),\\
D_\lambda&=\frac{\lambda}{(\lambda+1)^3}+O\left(e^{-\lambda}\right).
\end{align*}
Moreover, for each of these three constants, we were able to obtain bounds, which are very tight for large values of $\lambda$. For small values of $\lambda$, we described a method for estimating $C_\lambda$ and $B_\lambda$ with the help of a functional  equation for $M'_{\lambda}(x)$, which is simpler than the one for $M_\lambda(x)$ and satisfies \eqref{eq:9}. This remarkable property was then used to show that $M'_\lambda(x)\to M'(x)$ uniformly for all $x\in[0,\infty)$ as $\lambda\to 0$, which implies that $C_\lambda\to C$ and $B_\lambda\to B$ as $\lambda\to 0$. The constants $C$ and $B$ were obtained by R\'{e}nyi  in 1958 for the parking problem in the case of the uniform distribution.

{\bf Acknowledgements.}
The authors are thankful to Jack Hymowitz for the visualization of these results and to the Reviewer who suggested a possible improvement to the decay rate of the error terms in \eqref{eq:6} and \eqref{40}.

%
%



\end{document}